\documentclass[a4paper,10pt]{amsart}
\usepackage[arrow,matrix]{xy}
\usepackage{amsmath,amssymb,amscd,bbm,amsthm,mathrsfs,dsfont}
\theoremstyle{plain}
\newtheorem{theorem}{Theorem}[section]
\newtheorem{corollary}[theorem]{Corollary}
\newtheorem{lemma}[theorem]{Lemma}
\newtheorem{prop}[theorem]{Proposition}
\theoremstyle{definition}
\newtheorem{defn}[theorem]{Definition}
\newtheorem{exa}[theorem]{Example}

\theoremstyle{remark}

\newtheorem{remark}{Remark}[section]

\DeclareMathOperator{\Ext}{Ext} \DeclareMathOperator{\Hom}{Hom}

\begin{document}
\title[$A$-infinity structures related to bi-Koszul algebras]
{\bf $A$-infinity structures related to bi-Koszul algebras}


\author{Jun-Ru Si}
\address{(S\,i) Department of Mathematics, Zhejiang University, Hangzhou
310027, China} \email{sijunru@126.com}
\author{Di-Ming Lu}
\address{(Lu) Department of Mathematics, Zhejiang University, Hangzhou
310027, China}
\email{dmlu@zju.edu.cn}

\thanks {The work was supported by the NSFC (Grant No. 10571152) and
partially by
the NSF of Zhejiang Province of China (Grant No. J20080154)}
\date{}
\maketitle

\begin{abstract}
Let $A$ be a bi-Koszul algebra, we describe all possible
$A_\infty$-algebra structures on the Ext-algebra $E(A)$, and prove
that $E(A)$ must be $[m_2, m_3]$-finitely generated. An equivalent
description for a connected graded algebra to be a bi-Koszul algebra
is given in terms of $A_\infty$-language. The case that $E(A)$ is
endowed with minimal number of multiplications is discussed for
decomposition.
\end{abstract}

\vskip7mm
\section*{Introduction}

To understand certain homological properties of graded algebras
whose trivial modules admit non-pure resolutions, the authors
introduced what they have called bi-Koszul algebras in \cite{LS}.
Any non-Koszul Artin-Schelter regular algebras generated in degree 1
of global dimension four are the examples. Different from algebras
with certain pure resolutions of the trivial modules (such as Koszul
algebras \cite{P}, $d$-Koszul algebras \cite{Be}, piecewise-Koszul
algebras \cite{LHL}, \textit{etc}.) and $\mathcal{K}_2$-algebras
\cite{CS}, bi-Koszul algebras lose a nice homological property that
their Ext-algebras are finitely generated.

This may be remodeled if one endows ``generating'' with an
appropriate meaning. For example, Keller claimed that  Ext-algebras
are $A_\infty$-generated by their homogenous components of degree 1
for a large number of graded algebras \cite[Proposition 1(b)]{K1}.
Though it is nice to have finite generating components on one hand,
it maybe require, as a redeem, infinite multiplications to guarantee
finitely generating on the other hand. One of goals of this paper is
to find out the multiplications on a given Ext-algebra as less as
possible to carry out the finitely generating. We prove that the
Ext-algebra of any bi-Koszul algebra is $[m_2, m_3]$-finitely
generated.

We examine all possible $A_\infty$-algebra structures, corresponding
to a bi-Koszul algebra $A$ determined by $\varDelta_d$, to get at
most five multiplications $m_2,m_3,m_4,m_d,m_{d+1}$. An
$A_\infty$-version duality theory for a bi-Koszul algebra is given.
The case that $E(A)$ is endowed with minimal number of
multiplications $m_2,m_d,m_{d+1}$ is especially interesting for
decomposition. Two single $A_\infty$-algebras are obtained here and
can be returned to the $A_\infty$-structure of $E(A)$ by a bridge.

We introduce a modified concept of ``generating'' which reflects
some balance between multiplications and elements in the
$A_\infty$-algebra system and prove that there exists an
$A_\infty$-algebra structure on $E(A)$ of a bi-Koszul algebra $A$
such that $E(A)$ is $[m_2,m_3]$-finitely generated by $E^1(A)$,
$E^2(A)$ and $E^3(A)$. A new criterion for a bi-Koszul algebra to be
strongly is given in the $A_\infty$-version. Based on the fact that
the $A_\infty$-Ext-algebra $E(A)$ is unique up to quasi-isomorphism,
we discuss whether an $A_\infty$-algebra $E(A)$ is generated by
$E^1(A)$.

\vskip5mm
\section{$A_\infty$-algebras and bi-Koszul algebras}

In this section, we review basic material necessary for the paper:
Ext-algebras, $A_\infty$-algebras, and bi-Koszul algebras.

\subsection{Ext-algebras}

Throughout we fix a field $\mathds{F}$. We always assume that a
graded algebra $A=\mathds{F}\oplus A_1\oplus A_2\oplus\cdots$ is
locally finite, connected, and generated in degree 1. The graded
Jacobson radical of $A$, denoted by $J$, is $J=A_{\ge 1}$. Let
$Gr(A)$ denote the category of graded left $A$-modules. The
morphisms in this category, denoted by $\Hom_{Gr(A)}(M, N)$ for $M,
N\in Gr(A)$, are graded $A$-module maps of degree zero. For $M\in
Gr(A)$, we denote the $n^\mathrm{th}$ {\it shift\/} of $M$ by $M[n]$
where $M[n]_j=M_{j+n}$.

We write $\underline{\Ext}^*_A$ the derived functor of the graded
$\underline{\Hom}^*_A$ functor
$$
\underline{\Hom}^*_A(M, N):=\bigoplus_n \Hom_{Gr(A)}(M, N[n]),
$$
and denote
$$
E(A): = \underline{\Ext}_A^*(\mathds{F}, \mathds{F}),\quad  \quad
\mathsf{E}(M): = \underline{\Ext}_A^*(M, \mathds{F}),
$$
the {\it Koszul dual\/} of the algebra $A$ and the {\it Koszul
dual\/} of the module $M\in Gr(A)$, respectively. $E(A)$ is equipped
with a bigraded algebra structure by the Yoneda product with the
$(i, j)^\textmd{th}$ component $\underline{\Ext}_A^i(\mathds{F},
\mathds{F})_{-j}$, we also call it (classical) {\it Ext-algebra\/}
of $A$. Here, $i$ is the {\it cohomology degree} and $-j$ is the
{\it internal degree}. Note that the internal degree in $E(A)$ is
non-positive. For simplicity, we promise
$E^i_j(A):=\underline{\Ext}_A^i(\mathds{F}, \mathds{F})_{-j}$.
Similarly, $\mathsf{E}(M)$ is a bigraded left $E(A)$-module with the
$(i, j)^\textmd{th}$ component
$\mathsf{E}^i_j(M):=\underline{\Ext}_A^i(M, \mathds{F})_{-j}$.

The (classical) Ext-algebra $E(A)$ carries rich information about
the algebra $A$ and its module category, but it does not contain
enough information to recover the original algebra in general, the
``hidden'' information is revealed in the $A_\infty$-world.

\subsection{$A_\infty$-algebras}

There are different methods to give the definition of an
$A_\infty$-algebra (algebraical, geometrical, operadic, etc.), but
here we prefer the algebraical definition of an $A_\infty$-algebra.
We refer to \cite{K2} or \cite{LPWZ2} for the details.

\begin{defn}
An {\it $A_\infty$-algebra\/} over a field $\mathds{F}$ is a
$\mathbb{Z}$-graded vector space
 $$
 E=\bigoplus_{p\in \mathbb{Z}}E^p
 $$
endowed with a family of graded $\mathds{F}$-linear maps
 $$
 m_n: E^{\otimes n}\to E,\;\; (n\geq 1)
 $$
of degree $2-n$ satisfying the {\it Stasheff's identities:\/} for
all $n\geq 1$,

\medskip  \noindent \textsf{SI(n)} $\qquad\qquad \qquad\qquad$ $\sum
(-1)^{i+jt} m_{i+1+j}(1^{\otimes i}\otimes m_t \otimes 1^{\otimes
j})=0,$

\medskip \noindent where the sum runs over all decompositions $n=
i+t+j$ ($i$, $j\geq 0$ and $t\geq 1$).
\end{defn}

The graded maps $m_n$ for $n\ge 3$ are called {\it higher
multiplications\/} of $E$. An $A_\infty$-algebra $E$ is {\it
strictly unital\/} if $E$ contains an element 1 which acts as a
two-sided identity with respect to $m_2$, and for $n\ne 2,
m_n(x_1\otimes \cdots\otimes x_n)=0$ if $x_i=1$ for some $i$. The
$A_\infty$-algebras in this paper are always assumed to be strictly
unital. An $A_\infty$-algebra with zero $m_1$ is called {\it
minimal\/}. An {\it $A_\infty$-subalgebra\/} of $E$ is a graded
subspace $F$ such that $m_n$ maps $F^{\otimes n}$ to $F$ for all
$n\ge 1$. By an $A_{\infty}$-algebra $E$ being {\it generated\/} by
$E^1$ we mean that for any $p\geq2$,
$$
E^p=\sum m_l(E^{i_1}\otimes \cdots \otimes E^{i_l}),
$$
\noindent where the sum runs over all decompositions
$i_1+\cdots+i_l+2-l=p\;$ ($i_1, \cdots, i_l\geq1)$ and $l\geq1$.

Let $E$ and $F$ be two $A_{\infty}$-algebras. A {\it  morphism\/} of
$A_{\infty}$-algebras $f: E\to F$ is a family of graded
$\mathds{F}$-linear maps
$$f_n: E^{\otimes n}\to F, \;\; n\ge 1$$ of degree $1-n$ satisfying
the {\it Stasheff's morphism identities:\/} for all $n\geq 1$,

\medskip \noindent \textsf{  MI(n)} $ \qquad$ $\sum (-1)^{i+jt}
f_{i+1+j}(1^{\otimes i}\otimes m_t \otimes 1^{\otimes j})=\sum
(-1)^{w}m_r(f_{i_1}\otimes\cdots\otimes f_{i_r}),$

\medskip \noindent where the first sum runs over all decompositions
$n= i+t+j$ ($i$, $j\geq 0$ and $t\geq 1$), and the second sum runs
over all $1\leq r\leq n$ and all decompositions $n=i_1+\cdots+i_r$
(all $i_j\geq1 $); the sign on the right-hand side is given by
$w=(r-1)(i_1-1)+(r-2)(i_2-1)+\cdots  +(i_{r-1}-1)$.

An $A_\infty$-morphism in this paper is also required to be strictly
unital (see \cite{LPWZ1}). A morphism $f$ is called a {\it
quasi-isomorphism\/} if $f_1$ is a quasi-isomorphism. A morphism $f$
is called a {\it strict isomorphism\/} if $f_i=0$ for any $i\geq2$
and $f_1$ is an isomorphism.

$A_\infty$-algebras have been in use in topology since their
introduction by Stasheff. Their applicability in an algebraic
context was made clear by the minimality theorem, proven by
Kadeishvili \cite{Ka, K2}.

\begin{theorem}(The minimality theorem)
Let $E$ be an $A_\infty$-algebra. Then the cohomology $H^*E$ has an
$A_\infty$-algebra structure such that $m_1=0$, $m_2$ is induced by
$m^E_2$, and $H^*E$ is quasi-isomorphic to $E$ as
$A_{\infty}$-algebras. \qed
\end{theorem}

The techniques used to prove the minimality theorem all yield
explicit methods to compute an $A_\infty$-algebra structure. The
Ext-algebra $\underline{\Ext}^*_A(\mathds{F}, \mathds{F})$ is the
cohomology of $\mbox{End}_A(P)$, where $P$ is any free resolution of
$_A\mathds{F}$. Since $E=\mbox{End}_A(P)$ is a differential graded
algebra, by the minimality theorem,
$\underline{\Ext}_A^*(\mathds{F}, \mathds{F})$ has a natural
$A_\infty$-structure, which is called an {\it
$A_\infty$-Ext-algebra} of $A$.  By abuse of notation we still use
$E(A)$ to denote an $A_\infty$-Ext-algebra. The importance is that
the information from the $A_\infty$-algebra $E(A)$ is sufficient to
recover $A$.

We are mainly, in this paper, interested in the $A_\infty$-algebra
$E(A):=\bigoplus_{p, i\in \mathbb{Z}}E^p_i(A)$ which is bigraded
with the lower grading inherited from the graded algebra $A$. Each
multiplication $m_n$, as well as each morphism between two bigraded
$A_{\infty}$-algebras, must preserve the lower grading.

An $A_\infty$-algebra $E(A)$ that we consider in the paper always
comes from a free resolution. Different choice of the free
resolutions yields quasi-isomorphic $A_\infty$-algebra structures on
$E(A)$. Under the assumption on $A$, any choice of such an
$A_\infty$-algebra structure on $E(A)$ with the multiplications
$\{m_n\}_{n\geq1}$ has the following properties: $m_1=0$, $m_2$ is
the Yoneda product of $E(A)$, and $E^{2}(A)$ is $A_\infty$-generated
by $E^{1}(A)$; that is, $E^{2}_n(A)=m_n(
E^{1}(A)\otimes\cdots\otimes E^{1}(A) )$ for each $n\ge 2$.
Moreover, there exists an $A_\infty$-algebra structure on $E(A)$
such that $E(A)$ is generated by $E^1(A)$ (\cite{K1}). For more
properties we refer to \cite{K1,K2} or \cite{LPWZ1,LPWZ2}.

\subsection{Bi-Koszul algebras}

To extend Koszulity to a graded algebra with a bi-degree resolution
of the ground field, the authors introduced what they have called
bi-Koszul algebras in \cite{LS}.

\begin{defn}
A {\it bi-Koszul algebra\/} (determined by $\varDelta_d$) is a
connected graded algebra $A$ whose trivial module $\mathds{F}$ has a
minimal graded free resolution $\mathcal{P}$ such that each $P_n$ is
generated in degrees $\varDelta_d(n)$ for all $n\geq 0$, where the
degree distribution $ \varDelta_d: \mathbb{N} \rightarrow
\mathbb{N}\times \mathbb{N}$ is defined, for a fixed integer $d\ge
2$, by
$$
\varDelta_d(n)=\left\{\begin{array}{llll}
\frac{n}{3}(2d, 2d), & \mbox{if $n \equiv 0 (\textrm{mod} 3)$;}\\
\frac{n-1}{3}(2d, 2d)+(1,1), & \mbox{if $n \equiv 1 (\textrm{mod} 3)$;}\\
\frac{n-2}{3}(2d, 2d)+(d, d+1), & \mbox{if $n \equiv 2 (\textrm{mod}
3)$.}\end{array} \right.
$$
\end{defn}

For simplicity, $\varDelta_d(n)$ is used to express both of its
image $(x, y)$ and of the set $\{x, y\}$. Artin-Schelter regular
algebras of global dimension 4 of types (\emph{13431}) and
(\emph{12221}) are the examples by taking $d=2$ and $d=3$,
respectively. We refer to \cite{LS} for the details.


\begin{theorem} \label{thm1} \cite{LS}
The following statements are equivalent:
\begin{enumerate}
\item $A$ is a bi-Koszul algebra determined by $\varDelta_d$;
\item $E(A)$ begins with $E^1(A)=
E^1_1(A)$, $E^2(A)= E^2_d(A)\oplus E^2_{d+1}(A)$, $E^3(A)=
E^3_{2d}(A)$, and for each $n\geq 1$,
\begin{enumerate}
\item $E^{3n}(A)= \overbrace{E^3(A)E^3(A)\cdots E^3(A)}^n$,
\item $E^{3n+1}(A)= E^1(A)E^{3n}(A)=E^{3n}(A)E^1(A)$,
\item $E^{3n+2}(A)\cong E^2(A)E^{3n}(A)\oplus \mathsf{E}^2_{2nd+d+1}
(J\Omega^{3n}(\mathds{F}))$ as $\mathds{F}$-spaces.
\hfill{$\square$}
\end{enumerate}
\end{enumerate}
\end{theorem}

In the above theorem, the obstruction
$\mathsf{E}^2_{2nd+d+1}(J\Omega^{3n}(\mathds{F}))$ arises from the
bigger degree in $\varDelta_d(3n+2)$. We call a bi-Koszul algebra
$A$ {\it strongly\/} if the obstruction is vanished. In graded
algebras setting, it is clear that the Ext-algebra $E(A)$ of a
strongly bi-Koszul algebra $A$ is generated by $E^1(A), E^2(A)$ and
$E^3(A)$, but it is not sure that the Ext-algebra $E(A)$ of a
bi-Koszul algebra $A$ is finitely generated. There is a remedy of
finitely generating on $E(A)$ by using higher multiplications in
Section 3.


\vskip5mm
\section{$A_\infty$-Ext-algebras of bi-Koszul algebras}

In this section, we examine the possible multiplications on $E(A)$
as an $A_\infty$-algebra for a bi-Koszul algebra $A$ by using
information about the grading of $E(A)$. An $A_\infty$-version
duality theory of bi-Koszul algebras is given. In particular, we
discuss a kind of bi-Koszul algebras whose Ext-algebras are endowed
with the minimal number of nonzero multiplications.

\subsection{$A_\infty$-structures on $\mathbf{E(A)}$}

For the sake of convenience, we write
$$m_l(E^{t_1}\cdots
E^{t_l}):=m_l(E^{t_1}\otimes\cdots \otimes E^{t_l}).
$$

The following lemma gives an equivalent definition of the bi-Koszul
algebra which is characterized by its Ext-algebra.

\begin{lemma}\label{le0}
$A$ is a bi-Koszul algebra if and only if for any $n\geq0$,
$E^n_j(A)=0$ for $j\notin \varDelta_d(n)$.
\end{lemma}

\begin{proof}
Similar to the proof in \cite[Proposition 2.1.3]{BGS}.
\end{proof}

Before determining all possible multiplications on $E(A)$, we claim
that $m_2, m_d$ and $m_{d+1}$ must be non-trivial.

\begin{prop}\label{pp0}
Let $A$ be a bi-Koszul algebra determined by $\varDelta_d$. An
$A_\infty$-algebra $E:=E(A)$ must have nonzero multiplications
$m_2$, $m_d$ and $m_{d+1}$.
\end{prop}

\begin{proof}
As mentioned in the last section, $m_2$ is the Yoneda product, so we
need only to show that both $m_d$ and $m_{d+1}$ are non-trivial.
Noting that $E^2=E^2_d\oplus E^2_{d+1}$ and $E^2$ is generated by
$E^1$, we have
$$
m_d(E^{1}\cdots E^{1})=E^{2}_d,\;\;\; m_{d+1}(E^{1}\cdots
E^{1})=E^{2}_{d+1}.
$$
So we get $m_2$, $m_d$ and $m_{d+1}$ are nonzero.
\end{proof}

One of main results of this section is

\begin{theorem}\label{thm0}Let $A$ be a bi-Koszul algebra determined
by $\varDelta_d$. Then all possible non-trivial multiplications on
the $A_\infty$-Ext-algebra $E(A)$ are $m_2$, $m_3$, $m_4$, $m_d$ and
$m_{d+1}$.
\end{theorem}

\begin{proof}
Denote $E:=E(A)$. Let $m_l$ be a multiplication on $E(A)$. Since
only information about the grading of $E(A)$ is considered in the
following, we can neglect the order of $E^{i_1}, \cdots, E^{i_l}$
acted by $m_l$. Write
$$
M:=m_l(E^{3k_1+t_1}\cdots E^{3k_\alpha+t_\alpha}
E^{3k_{\alpha+1}+2}\cdots E^{3k_l+2})
$$
where $\alpha\leq l$ and $t_j=0$ or $1$ $(1\leq j \leq\alpha)$.
Denote $\beta=l-\alpha$. So
$$
M\subseteq E^{3(k_1+\cdots +k_l)+(t_1+\cdots+
t_\alpha)+2\beta+2-l}
$$
and the lower grading of $M$ falls into the set
$$
\{2d(k_1+\cdots +k_l)+(t_1+\cdots+ t_\alpha)+d \beta +j\;\;|\;\;
j=0,1,\cdots, \beta\},
$$
where $0\leq t_1+\cdots +t_\alpha\leq l-\beta$.

(1) If $(t_1+\cdots t_\alpha)+2\beta+2-l=3k$ $(k\geq0)$, then
$$
E^{3(k_1+\cdots+k_l)+(t_1+\cdots+t_\alpha)+2\beta+2-l}=E^{3(k_1+\cdots
+k_l)+3k}_{2d(k_1+\cdots+ k_l)+2dk}.
$$
We have the following inequalities:
$$
(t_1+\cdots+t_\alpha)+d\beta\leq 2kd \leq (t_1+\cdots
+t_\alpha)+d\beta+\beta,
$$
which produce the solutions of $(k, \beta, l)$ as the following list
$$
(0, 0, 2), \;(1, 1, d),\; (1, 1, d+1),\; (1, 2, 3),\; (2, 4,
4).\eqno(S1)
$$

(2) If $(t_1+\cdots t_\alpha)+2\beta+2-l=3k+1$ $(k\geq0)$, then
$$
E^{3(k_1+\cdots+k_l)+(t_1+\cdots+t_\alpha)+2\beta+2-l}=E^{3(k_1+\cdots
+k_l)+3k+1}_{2d(k_1+\cdots+ k_l)+2dk+1}.
$$
The inequalities:
$$
(t_1+\cdots +t_\alpha)+d\beta\leq 2kd+1 \leq (t_1+\cdots +t_\alpha)
+d\beta+\beta
$$
imply the solutions of $(k, \beta, l)$ in the following list
$$ (0, 0, 2),\;(1, 2, 2), \;(1, 2, 3).\eqno(S2)$$

(3) If $(t_1+\cdots t_\alpha)+2\beta+2-l=3k+2$ $(k\geq0)$, then
\begin{eqnarray*}
&&E^{3(k_1+\cdots +k_l)+(t_1+\cdots+t_\alpha)+2\beta+2-l}
\\&=&E^{3(k_1+\cdots +k_l)+3k+2}_{2d(k_1+\cdots+ k_l)+2dk+d}\oplus
E^{3(k_1+\cdots +k_l)+3k+2}_{2d(k_1+\cdots+
k_l)+2dk+d+1}.
\end{eqnarray*}
We have the following inequalities:
$$
(t_1+\cdots +t_\alpha)+d\beta\leq 2kd+d \leq (t_1+\cdots
+t_\alpha)+d\beta+\beta,
$$
or
$$
(t_1+\cdots +t_\alpha)+d\beta\leq 2kd+d+1 \leq (t_1+\cdots +t_\alpha)
+d\beta+\beta.
$$
The solutions of $(k, \beta, l)$ are listed in the following
$$(0, 1, 2),\; (1, 3, 3),\;(0, 0, d),\; (0, 1, 2),\; (0, 1, 3),
\; (0, 0, d\!+\!1),\; (1, 3, 3),\; (1, 3, 4).\eqno(S3)$$

In conclusion of (S1)-(S3), all possible solutions of $l$ are $2, 3,
4, d$ and $d+1$. This completes the proof.
\end{proof}

\begin{corollary}
Let $A$ be a bi-Koszul algebra determined by $\varDelta_d$.
\begin{enumerate}
\item If $d=2$ or $3$, the possible non-trivial
multiplications on $E(A)$ are $m_2$, $m_3$ and $m_4$.
\item If $d=4$, the possible non-trivial
multiplications on $E(A)$ are $m_2$, $m_3$, $m_4$ and $m_5$.
\item If $d\geq5$, the possible non-trivial
multiplications on $E(A)$ are $m_2$, $m_3$, $m_4$, $m_d$ and
$m_{d+1}$.  \hfill{$\square$}
\end{enumerate}
\end{corollary}

To describe what components the multiplications act on non-trivial,
we denote
 $$
 E^{[0]}:=\bigoplus_{k\geq0} E^{3k},\;\;
 E^{[1]}:=\bigoplus_{k\geq0}E^{3k+1},\;\;
 E^{[2]}:=\bigoplus_{k\geq0}E^{3k+2}=E^{[2]}_{(d)}\oplus E^{[2]}_{(d+1)}.
 $$

The following proposition is clear from the proof of Theorem
\ref{thm0}.

\begin{prop}\label{pp0}Let $A$ be a bi-Koszul algebra determined
by $\varDelta_d$, $E$ the Ext-algebra of $A$. Then the possible
nonzero components of $m_i$ $(i=2, 3, 4, d, d+1)$ are:

\medskip{\small
\begin{center}
\begin{tabular}{c|ccc}
$*$& fall into $E^{[0]}$ &  fall into $E^{[1]}$ & fall into
$E^{[2]}$
\\
\hline $m_2$ & $E^{[0]}E^{[0]}$ & $E^{[0]}E^{[1]}$,
$E^{[2]}_{(d)}E^{[2]}_{(d+1)}$ & $E^{[0]}E^{[2]}$ \\
$m_3$ & $E^{[0]}E^{[2]}_{(d)}E^{[2]}_{(d)}$ &
$E^{[1]}E^{[2]}_{(d)}E^{[2]}_{(d)}$ & $E^{[0]}E^{[1]}E^{[2]}_{(d)}$,
$E^{[2]}E^{[2]}_{(d)}E^{[2]}_{(d)}$ \\
$m_4$ & $E^{[2]}_{(d)}E^{[2]}_{(d)}E^{[2]}_{(d)}E^{[2]}_{(d)}$ &  &
$E^{[1]}E^{[2]}_{(d)}E^{[2]}_{(d)}E^{[2]}_{(d)}$ \\
$m_d$ & $\underbrace{E^{[1]}\cdots
E^{[1]}}_{d-1}E^{[2]}_{(d+1)}$ &  & $\underbrace{E^{[1]}\cdots
E^{[1]}}_d$ \\
$m_{d\!+\!1}$ & $\underbrace{E^{[1]}\cdots
E^{[1]}}_{d}E^{[2]}_{(d)}$ &
& $\underbrace{E^{[1]}\cdots E^{[1]}}_{d+1}$ \\
\end{tabular}\end{center}

\hskip8mm $*:$\;  including all permutations of the components
listed above.} \qed
\end{prop}

\begin{defn}
We call an $A_\infty$-algebra $E=(E; m_2, m_3, m_4, m_d, m_{d+1})$
{\it reduced\/}, if all possible nonzero components of
multiplications are in the above table.
\end{defn}

\begin{corollary}\label{cc1}
$A$ is a bi-Koszul algebra  if and only if any $A_\infty$-algebra
structure on $E(A)$ is reduced.
\end{corollary}

\begin{proof} The necessity is from Theorem \ref{thm0}
and Proposition \ref{pp0}. Now suppose that any $A_\infty$-algebra
on $E(A)$ is reduced. Take the $A_\infty$-algebra $E(A)$ that is
generated by $E^1(A)$, then $A$ is a bi-Koszul algebra by checking
the lower grading.
\end{proof}

\subsection{Truncated bi-Koszul algebras}

We discuss a kind of bi-Koszul algebras whose Ext-algebras are
endowed with the minimal number of non-trivial multiplications $m_2,
m_d$ and $m_{d+1}$.

\begin{defn}\label{def5}
Let $A$ be a bi-Koszul algebra determined by $\varDelta_d$,
$E:=E(A)$ its Ext-algebra. We say that $A$ is {\it truncated\/} if
the $A_\infty$-Ext-algebra $E$ only has the non-trivial
multiplications $m_2, m_d, m_{d+1}$ and the possible nonzero actions
of $m_i$ $(i=2, d, d+1)$ are on
\medskip{\small
\begin{center}
\begin{tabular}{c|ccc}
$*$& fall into $E^{[0]}$ &  fall into $E^{[1]}$ & fall into
$E^{[2]}$
\\
\hline $m_2$ & $E^{[0]}E^{[0]}$ & $E^{[0]}E^{[1]}$,
$E^{[2]}_{(d)}E^{[2]}_{(d+1)}$ & $E^{[0]}E^{[2]}$ \\
$m_d$ & $\underbrace{E^{[1]}\cdots E^{[1]}}_{d-1}E^{[2]}_{(d+1)}$ &
& $\underbrace{E^{[1]}\cdots
E^{[1]}}_d$ \\
$m_{d\!+\!1}$ & $\underbrace{E^{[1]}\cdots
E^{[1]}}_{d}E^{[2]}_{(d)}$ &
& $\underbrace{E^{[1]}\cdots E^{[1]}}_{d+1}$ \\
\end{tabular}\end{center}

\hskip18mm $*:$\;  including all permutations of the components
listed above.}

\noindent By abuse of notation we also say that the
$A_\infty$-algebra $(E; m_2, m_d, m_{d+1})$ is {\it truncated\/} in
this case.
\end{defn}

\begin{exa}
All Artin-Schelter regular algebras listed in
\cite[Theorem A]{LPWZ2} are truncated bi-Koszul algebras.
\end{exa}

The following result comes from Proposition \ref{pp1} of Section 3.

\begin{prop}\label{tp1}
A truncated bi-Koszul algebra must be strongly.
\end{prop}

A minimal $A_\infty$-algebra is called {\it single\/} if it has only
one non-trivial higher multiplication (\textit{i.e.} a $(2,
p)$-algebra discussed in \cite{HL}). Single $A_\infty$-algebras are
related to $p$-Koszul algebras (\cite{HL}).

A bigraded algebra is called {\it pure\/} in the sense that every
component is supported in a single lower grading. We also say that a
minimal $A_\infty$-algebra which is bigraded is {\it pure\/} if the
underlying bigraded algebra itself is pure.

We need the following lemma.  Consider an $A_\infty$-algebra $(E;
m_2, m_d, m_{t})$ $(2<d< t$ and $2+t\neq 2d$). All non-trivial
Stasheff's identities, in this case, are listed as follows:

\begin{itemize}
{\small
\item [\textsf{SI(3)}:] $m_2(m_2\otimes 1)=m_2(1\otimes m_2)$;
\item [\textsf{SI(d+1)}:]$\sum_{i+j=d-1}(-1)^{i}m_d(1^{\otimes i}\otimes m_2\otimes
1^{\otimes j})=m_2(1\otimes m_d)-(-1)^{d}m_2(m_d\otimes 1)$;
\item [\textsf{SI(t+1)}:]$\sum_{i+j=t-1}(-1)^{i}m_t(1^{\otimes i}\otimes m_2\otimes
1^{\otimes j})=m_2(1\otimes m_t)-(-1)^{t}m_2(m_t\otimes 1)$;
\item [\textsf{SI(2d-1)}:]$\sum_{i+j=d-1}(-1)^{i+dj}m_d(1^{\otimes i}\otimes m_d\otimes
1^{\otimes j})=0$;
\item [\textsf{SI(d+t-1)}:]$\sum_{i+j=d-1}(-1)^{i+tj}m_d(1^{\otimes i}\otimes m_t\otimes
1^{\otimes j})=\sum_{i+j=t-1}(-1)^{i+dj+1}m_t(1^{\otimes i}\otimes
m_d\otimes 1^{\otimes j})$;
\item [\textsf{SI(2t-1)}:]$\sum_{i+j=t-1}(-1)^{i+tj}m_t(1^{\otimes i}\otimes m_t\otimes
1^{\otimes j})=0$.}
\end{itemize}

\begin{lemma}\label{tl1}
Let $E$ be a connected graded algebra with three graded
$\mathds{F}$-linear maps $m_n: E^{\otimes n}\to E$ $(n=2, d, t)$.
Suppose $2<d<t$ and $2+t\neq 2d$. Then the following statements are
equivalent.
\begin{enumerate}
\item $(E; m_2, m_d, m_t)$ is an $A_\infty$-algebra;
\item $E$ together with $\{m_2, m_d, m_t\}$ satisfies
\begin{itemize}{ \item [(a)] $(E; m_2, m_d)$ is single,
\item [(b)] $(E; m_2, m_t)$ is single, \item [(c)] $m_d$ and $m_t$
obey \textsf{SI(d+t-1)}.}
\end{itemize}
\end{enumerate}
\end{lemma}

\begin{proof} It is clear by noting that all non-trivial Stasheff's
identities of a single
$A_\infty$-algebra with the higher multiplication $m_p$ are
$\textsf{SI(3)}, \textsf{SI(p+1)}, \textsf{SI(2p-1)}$.
\end{proof}

From the lemma above, one may decompose a truncated
$A_\infty$-algebra into two single $A_\infty$-algebras, while the
Stasheff's identity \textsf{SI(2d)} serves as a bridge between two
single $A_\infty$-algebras.

\begin{prop}\label{tp2}
Let $A$ be a truncated bi-Koszul algebra determined by $\varDelta_d$
($d\ge 4$). Then both $(E; m_2, m_d)$ and $(E; m_2, m_{d+1})$ are
single $A_\infty$-algebras,  generated by $E^{1}$, $E^{2}$ and
$E^{3}$.
\end{prop}

\begin{proof}
This is a direct result of Proposition \ref{tp1} and Lemma
\ref{tl1}.
\end{proof}

The single $A_\infty$-algebra $(E; m_2, m_d)$ in the proposition
above is not the $A_\infty$-Ext-algebra of any graded algebra, since
there are no components acted by $m_2$ and $m_d$ that fall into
$E^{2}_{d+1}$; neither is $(E; m_2, m_{d+1})$.

Drawing upon the $A_\infty$-algebras $(E; m_2, m_d)$ and $(E; m_2,
m_{d+1})$, the $A_\infty$-Ext-algebra $E$ can be decomposed further
into two single $A_\infty$-algebras which are both pure.

\begin{theorem}\label{tt1}
Let $A$ be a truncated bi-Koszul algebra determined by $\varDelta_d$
($d\ge 4$), $E:=E(A)$ its Ext-algebra. Set
$$
F:=E^{[0]}\oplus E^{[1]}\oplus E^{[2]}_{(d)} ,\;\;\;G:=E^{[0]}\oplus
E^{[1]}\oplus E^{[2]}_{(d+1)}.
$$
Then
\begin{enumerate}
\item $(F; m_2, m_d)$ is a pure and single $A_\infty$-subalgebra of $(E; m_2,
m_d)$, where $m_d$ is determined by $m_2$ and $m_d
\mid_{(E^1)^{\otimes d}}$;
\item [(2)] $(G; m_2, m_{d+1})$ is a pure and single $A_\infty$-subalgebra of
$(E; m_2, m_{d+1})$, where $m_{d+1}$ is determined by $m_2$ and
$m_{d+1}\mid_{(E^1)^{\otimes d+1}}$.
\end{enumerate}
\end{theorem}

\begin{proof}
It is easy to justify that both $(F; m_2, m_d)$ and $(G; m_2,
m_{d+1})$ are pure and single $A_\infty$-subalgebras of $(E; m_2,
m_d)$ and $(E; m_2, m_{d+1})$, respectively. And the nonzero actions
of $m_d$ are only on $(E^{[1]})^{\otimes d}$ and $m_{d+1}$ only on
$(E^{[1]})^{\otimes d+1}$.

Write $ m_2(x, y)$ by $xy$ or $x\cdot y $. Set $t=d$ or $d+1$. For
any homogeneous elements $x_1, \cdots, x_t\in E^{[1]}$ and $u\in
E^{[0]}$, we note that
\begin{eqnarray*}
&&m_t(x_1, \cdots, x_{i-1}, ux_i, \cdots, x_t)=m_t(x_1, \cdots,
x_{i-1}u, x_i, \cdots, x_t), (2\leq i\leq t),
\\&& m_t(ux_1, \cdots, x_t)=u\cdot
m_t(x_1,  \cdots, x_t),\\
&& m_t(x_1, \cdots, x_t\cdot u)=m_t(x_1, \cdots, x_t)\cdot u.
\end{eqnarray*}
Since $E^{3n+j}=E^{j}E^{3n}=E^{3n}E^{j}\;(j=1,2)$ by Proposition
\ref{tp1}, for any $x\in E^{3n+j}$ $(n\geq0)$, there exist $f,g\in
E^{j}$ and $u,v\in E^{3n}$ such that $x=fu=vg$. If $n=0$, set
$u=v=1$.

For any $x_1,\cdots,x_d\in E^{[1]}$, choose $z_1,\cdots,z_d\in E^1$
satisfying
 $$
 x_d=y_dz_d,x_{d-1}y_{d}=y_{d-1}z_{d-1},\cdots,
 x_{2}y_{3}=y_{2}z_{2},x_{1}y_{2}=y_{1}z_{1}
 $$
where $y_1,\cdots,y_d\in E^{[0]}$. We have
\begin{eqnarray*}
m_d(x_1, \cdots, x_d)&=&m_d(x_1, \cdots, x_{d-1}, y_dz_d)
\\&=&m_d(x_1, \cdots, x_{d-1}y_d, z_d)
\\&=&m_d(x_1, \cdots, y_{d-1}z_{d-1}, z_d)
\\&=&\cdots
\\&=&m_d(x_1y_{2},z_{2}, z_3\cdots, z_d)
\\&=&y_1\cdot m_d(z_1, \cdots, z_d),
\end{eqnarray*}
so $m_d$ is determined by $m_2$ and $m_d\mid_{(E^1)^{\otimes d}}$.

By the same method, $m_{d+1}$ is determined by $m_2$ and
$m_{d+1}\mid_{(E^1)^{\otimes d+1}}$. We complete the proof.
\end{proof}

In the $A_\infty$-Ext-algebra $(E;\{m_i\})$ of a graded algebra $A$,
$m_2$ is the Yoneda product and $m_i|_{(E^1)^{\otimes i}}$ can be
computed out explicitly for every $i\geq3$. This was demonstrated in
\cite{HL}. More concretely, the single higher multiplication in
either $(F; m_2, m_d)$ or $(G; m_2, m_{d+1})$ becomes definite in
form.

Now, we turn to find a way in which a truncated $A_\infty$-algebra
can be formed by jointing two single $A_\infty$-algebras together as
follows.

Suppose that $(E; m_2)$ is a bigraded algebra starting with
$E^{1}=E^{1}_{1}$, $E^{2}=E^{2}_{d}\oplus E^{2}_{d+1}$,
$E^{3}=E^{3}_{2d}$, and satisfying $E^{3n+i}=E^iE^{3n}=E^{3n}E^i$
for all $n\ge 1,\; i=1, 2, 3$. Define two single $A_\infty$-algebras
$(E; m_2, m_d)$ and $(E; m_2, m_{d+1})$ such that
\begin{itemize}
\item [(1)] the nonzero actions of $m_d$ are only on
$(E^{[1]})^{\otimes d}$ and $E^{[1]}\cdots E^{[2]}_{(d+1)}\cdots
E^{[1]}$ (including all permutations);
\item [(2)] the nonzero actions of $m_{d+1}$ are only on $(E^{[1]})^{\otimes
d+1}$ and $E^{[1]}\cdots E^{[2]}_{(d)}\cdots E^{[1]}$ (including all
permutations).
\end{itemize}
Then we have

\begin{prop}\label{tp3}
Let $(E; m_2)$ and $m_d, m_{d+1}$ be as above with $d\geq4$. If
$m_d$ is compatible with $m_{d+1}$ by \textsf{SI(2d)} on
$(E^{1})^{\otimes 2d}$. Then $(E; m_2, m_d, m_{d+1})$ is a truncated
$A_\infty$-algebra.
\end{prop}

\begin{proof}
We only need to show that $m_d$ is compatible with $m_{d+1}$ by
\textsf{SI(2d)} on $(E^{[1]})^{\otimes 2d}$ by Lemma \ref{tl1} and
the nonzero actions of $m_d$ and $m_{d+1}$.

Write $m_2(x, y)$ by $xy$ or $x\cdot y $. Note the nonzero actions
of $m_d$ and $m_{d+1}$ and the proof of Theorem \ref{tt1}. For any
$x_1,\cdots,x_{2d}\in E^{[1]}$, choose $z_1,\cdots,z_{2d}\in E^1$
satisfying
 $$
 x_{2d}=y_{2d}z_{2d},\;x_{2d-1}y_{2d}=y_{2d-1}z_{2d-1}, \;\cdots,\;
 x_{2}y_{3}=y_{2}z_{2},\;x_{1}y_{2} =y_{1}z_{1}
 $$
where $y_1,\cdots,y_{2d}\in E^{[0]}$.  We
have
\begin{eqnarray*}
&&m_d(x_1, \cdots, x_i, m_{d+1}(x_{i+1}, \cdots, x_{i+d+1}),
x_{i+d+2}, \cdots, x_{2d})
\\&=&m_d(x_1, \cdots, x_i, m_{d+1}(x_{i+1}, \cdots, x_{i+d+1}), x_{i+d+2},
\cdots, x_{2d-1}y_{2d}, z_{2d})
\\&=&m_d(x_1, \cdots, x_i, m_{d+1}(x_{i+1}, \cdots, x_{i+d+1}), x_{i+d+2}
y_{i+d+3}, \cdots, z_{2d})
\\&=&m_d(x_1, \cdots, x_i, m_{d+1}(x_{i+1}, \cdots, x_{i+d+1})\cdot
y_{i+d+2}, z_{i+d+2}, \cdots, z_{2d})
\\&=&m_d(x_1, \cdots, x_i, m_{d+1}(x_{i+1}, \cdots, x_{i+d+1}y_{i+d+2}),
z_{i+d+2}, \cdots, z_{2d})
\\&=&m_d(x_1, \cdots, x_i, m_{d+1}(x_{i+1}y_{i+2}, \cdots, z_{i+d+1}),
z_{i+d+2}, \cdots, z_{2d})
\\&=&m_d(x_1, \cdots, x_i, y_{i+1}\cdot m_{d+1}(z_{i+1}, \cdots,
z_{i+d+1}),
z_{i+d+2}, \cdots, z_{2d})
\\&=&m_d(x_1y_2, \cdots, z_i, m_{d+1}(z_{i+1}, \cdots, z_{i+d+1}),
z_{i+d+2},
\cdots, z_{2d})
\\&=&y_1\cdot m_d(z_1, \cdots, z_i, m_{d+1}(z_{i+1}, \cdots, z_{i+d+1}),
z_{i+d+2}, \cdots, z_{2d}).
\end{eqnarray*}
Similarly,
\begin{eqnarray*}
&&m_{d+1}(x_1, \cdots, x_i, m_{d}(x_{i+1}, \cdots, x_{i+d}),
x_{i+d+1}, \cdots, x_{2d})
\\&=&y_1\cdot m_{d+1}(z_1, \cdots, z_i, m_{d}(z_{i+1}, \cdots, z_{i+d}),
z_{i+d+1}, \cdots, z_{2d}).
\end{eqnarray*}

Set
 $$
 \varphi:=\sum_{i+j=d-1}(-1)^{i+(d+1)j}m_d(1^{\otimes i}\otimes m_{d+1}\otimes
 1^{\otimes j})+\sum_{i+j=d}(-1)^{i+dj}m_{d+1}(1^{\otimes i}\otimes
 m_{d}\otimes 1^{\otimes j}),
 $$
then $\varphi(x_1\otimes \cdots \otimes x_{2d}) =y_1\cdot
\varphi(z_1\otimes \cdots \otimes z_{2d})=0$ by the assumption.
Therefore, \textsf{SI(2d)} holds on all $(E^{[1]})^{\otimes 2d}$.

We complete the proof.
\end{proof}

\begin{remark}
In the proofs of Theorem \ref{tt1} and Proposition \ref{tp3}, we
ignore the $\sum$ when run up the multiplication $m_2$. This does
not affect the results.
\end{remark}

We finally give a condition under which the $A_\infty$-Ext-algebra
$E(A)$ is generated by $E^1(A)$.

\begin{prop}
Let $A$ be a truncated bi-Koszul algebra determined by
$\varDelta_d$ $(d\ge 4)$, $E:=E(A)$ its Ext-algebra. If either $(E;
m_2, m_d)$ or $(E; m_2, m_{d+1})$ is generated by $E^{1}$ and
$E^{2}$. Then $(E; m_2, m_d, m_{d+1})$ is generated by $E^{1}$.
\end{prop}

\begin{proof}
By Proposition \ref{tp1}, $E$ is generated by $E^1$, $E^2$, $E^3$ as
an associative algebra. Clearly, $E^2=m_d(\underbrace{E^1 \cdots
E^1}_d)+m_{d+1}(\underbrace{E^1\cdots E^1}_{d+1})$ in $(E; m_2, m_d,
m_{d+1})$. Therefore, to prove the result, we need only to show that
$E^3$ can be generated by $E^1$ and $E^2$ in $(E; m_2, m_d,
m_{d+1})$, which follows either from the assumption of $(E; m_2,
m_d)$ generated by $E^{1}$ and $E^{2}$ then
$$
E^{3}=\sum_{i=1}^d m_d(\underbrace{E^{1}\cdots E^{2}_{d+1}}_i\cdots
E^{1}),
$$
or from the assumption of $(E; m_2, m_{d+1})$ generated by $E^{1}$
and $E^{2}$ then
$$
E^{3}=\sum_{i=1}^{d+1} m_{d+1}(\underbrace{E^{1}\cdots
E^{2}_{d}}_i\cdots E^{1}).
$$
We get the result.
\end{proof}

\vskip5mm
\section{Balanced generating}

For a bi-Koszul algebra $A$, it is a question whether $E(A)$ is
finitely generated as a graded algebra. In this section, we
generalize the concept of ``generating'', and show that $E(A)$ is
$[m_2, m_3]$-finitely generated by $E^1(A), E^2(A)$ and $E^3(A)$ for
any bi-Koszul algebra $A$. An equivalent statement of a bi-Koszul
algebra is given in terms of such concept.

\subsection{$\mathbf{[m_2, m_3]}$-Generating}

The original concept of ``generating'' in the associative algebra
setting is defined with respect to the multiplication. When we work
in the field of $A_\infty$-algebras, we need a generalized concept
of ``generating'' to reflect certain balance between multiplications
and elements.

\begin{defn}
Let $E$ be an $A_\infty$-algebra. Suppose there exists a fixed
integer $l$ and multiplications $m_{n_1}, \cdots, m_{n_t} $ such
that, for any $p>l$,
$$
E^p=\sum_{\substack{k_1+\cdots+k_{n_i}\!+2-n_i=p\\ k_1, \cdots,
k_{n_i}\ge 1;\;\; 1\le i\le t}} m_{n_i}(E^{k_1}\otimes \cdots
\otimes E^{k_{n_i}}).
$$
We say that $E$ is {\it $[m_{n_1},\cdots,m_{n_t}]$-finitely
generated by $E^1, \cdots, E^l$.}
\end{defn}

\begin{remark}
In the case of $m_1=0$ and $[m_{n_1},\cdots,m_{n_t}]=[m_2]$, the
concept is the original one of finitely generating as an associative
graded algebra.
\end{remark}

Here are the examples.

If $A$ is a  $p$-Koszul algebra ($p\geq3$), then any
$A_\infty$-algebra $(E(A); m_2, m_p)$ is $[m_2, m_p]$-finitely
generated  by $E^1(A)$. This is obtained by noting the facts that
$E(A)$ is generated by $E^1(A)$ and $E^2(A)$, while
$E^2(A)=m_p(E^1(A)\otimes\cdots\otimes E^1(A))$ (\cite[Theorem
2.5]{HL}).

If $A$ is a bi-Koszul algebra, then there exists an
$A_\infty$-algebra $(E(A); m_2, m_3, m_4,\newline m_d, m_{d+1})$
such that $E(A)$ is $[m_2, m_3, m_4, m_d, m_{d+1}]$-finitely
generated by $E^1(A)$.

If we admit the set of multiplications to be infinite, Keller's
result tells that there exists an $A_\infty$-algebra structure on
$E(A)$ which is generated by $E^1(A)$ \cite[Proposition 1(b)]{K1}.

It is natural to expect that the multiplications in the set
$\{m_{n_1}, \cdots, m_{n_t}\}$ as less as possible. When $A$ is a
bi-Koszul algebra, though we can not claim $E(A)$ is
$[m_2]$-finitely generated, it does be finitely generated as long as
to add one higher multiplication.

Now let $A$ be a bi-Koszul algebra determined by $\varDelta_d$, and
$E:=E(A)$  the Ext-algebra of $A$ in the following. Denote

\begin{itemize}
\item $\mathbf{U}^{3k+2}$: the sum of the actions of $m_3$ on all
permutations of $E^{3k_1}$, $E^{3k_2+1}$ and $E^{3k_3+2}_{2dk_3+d}$
for any $k_1\geq1, k_2, k_3\geq0$ with $k_1+k_2+k_3=k$;
\item $\mathbf{V}^{3k+2}$: the sum of the actions of $m_3$ on all
permutations of
$E^{3k_1+2}_{2dk_1+d}$, $E^{3k_2+2}_{2dk_2+d}$ and
$E^{3k_3+2}_{2dk_3+d+1}$ for any $k_1, k_2, k_3\geq0$ with
$k_1+k_2+k_3=k-1$;
\item $\mathbf{W}^{3k+2}$:  the sum of the actions of $m_4$ on all
permutations  of
$E^{3k_1+1}$, $E^{3k_2+2}_{2dk_2+d}$, $E^{3k_3+2}_{2dk_3+d}$ and
$E^{3k_4+2}_{2dk_4+d}$ for any $k_1, k_2, k_3, k_4\geq0$ with
$k_1+k_2+k_3+k_4=k-1$.
\end{itemize}

\begin{prop}\label{pp2}Let  $A$  be a bi-Koszul algebra. Then for any
$k\geq1$,
$$
\begin{array}{ccclcrccc}
E^{3k+3}=m_2(E^{3}E^{3k})=m_2(E^{3k}E^{3});\\
E^{3k+1}=m_2(E^{1}E^{3k})=m_2(E^{3k}E^{1});\\
E^{3k+2}_{2dk+d}=m_2(E^{2}_{d}E^{3k})=m_2(E^{3k}E^{2}_d).
\end{array}
$$
\end{prop}

\begin{proof}
By Theorem \ref{thm1}, we need only to show the last equality:
$$m_2(E^{2}_{d}E^{3k})=m_2(E^{3k}E^{2}_d).$$

Since $m_2(E^{2}_{d}E^{3k})=m_2(m_d(E^{1}\cdots E^{1})E^{3k})$, to
get $m_2(E^{2}_{d}E^{3k})\subseteq m_2(E^{3k}E^{2}_d)$ we need only
to verify $m_2(m_d(x_1,\cdots,x_d),y)\in m_2(E^{3k}E^{2}_d)$ for any
$x_1, \cdots, x_d\in E^{1}$ and $y\in E^{3k}$. This is performed by
using the Stasheff's identity \textsf{SI(d+1)}
$$
m_2(m_d(x_1,\cdots,x_d),y)=m_d(x_1,\cdots,x_{d-1},m_2(x_d,y))
$$
with $m_2(x_d,y)\in E^{3k+1}$. Since $E^{3k+1}=m_2(E^{3k}E^{1})$,
\begin{eqnarray*}
m_d(x_1,\cdots,x_{d-1},m_2(x_d,y))&=&m_d(x_1,\cdots,x_{d-1},m_2(y',x'_d))
\\&=&m_d(x_1,\cdots,m_2(x_{d-1},y'),x'_d)
\end{eqnarray*}
with $y'\in E^{3k}$ and $x'_d\in E^{1}$. We can continue the
foregoing procedure to obtain
\begin{eqnarray*}
m_d(x_1, \cdots, x_{d-1}, m_2(x_d,y))&=&m_d(m_2(z,x'_1), \cdots,
x'_d)
\\&=&m_2( z,m_d(x'_1, \cdots, x'_d))
\end{eqnarray*}
with $z\in E^{3k}$ and $x'_1, \cdots, x'_d\in E^{1}$. Thus,
$m_2(E^{2}_{d}E^{3k})\subseteq m_2(E^{3k}E^{2}_d)$.

The converse $m_2(E^{3k}E^2_d)\subseteq
E^{3k+2}_{2dk+d}=m_2(E^2_dE^{3k})$ is clear from Theorem \ref{thm1}
again.
\end{proof}

\begin{lemma}\label{le1}Let $A$ be a bi-Koszul algebra.
Assume $k_i\geq0$ $(i=1,\cdots, d+1)$ and
$k_1+\cdots+k_{d+1}=k\geq1$. Then
$$
m_{d+1}(E^{3k_1+1}\cdots E^{3k_{d+1}+1}) \subseteq \sum_{\substack{
i_1+i_2=k, \\ i_1\geq1, i_2\geq0 }} m_2(E^{3i_1}E^{3i_2+2})+
\mathbf{U}^{3k+2}.
$$
\end{lemma}

\begin{proof}
There exists an integer $i$ $(1\leq i\leq d+1)$ such that
$k_i\geq1$. If $2\leq i \leq d+1$, using \textsf{SI(d+2)}, we get
\begin{eqnarray*}
&&m_{d+1}(E^{3k_1+1}\cdots E^{3k_{d+1}+1})\\
&=&m_{d+1}(E^{3k_1+1}\cdots m_2(E^{3k_i}E^{1})\cdots E^{3k_{d+1}+1})\\
&\subseteq& m_{d+1}(E^{3k_1+1}\cdots E^{3k_{i-1}+3k_i+1}E^{1}\cdots
E^{3k_{d+1}+1})+\mathbf{U}^{3k+2}\\
&\subseteq& m_{d+1}(E^{3(k_1+\cdots+k_i)+1}\cdots E^{1}E^{1}\cdots
E^{3k_{d+1}+1})+\mathbf{U}^{3k+2}
\end{eqnarray*}
with $ k_1+\cdots+ k_i\geq 1$.

So we can assume $k_1\geq1$. Using \textsf{SI(d+2)} again,
\begin{eqnarray*}
&&m_{d+1}(E^{3k_1+1}\cdots E^{3k_{d+1}+1})
\\&=&m_{d+1}(m_2(E^{3k_1}E^{1})E^{3k_2+1}\cdots E^{3k_{d+1}+1})
\\&\subseteq&
m_3(E^{3k_1}E^{3k_2+\cdots+3k_d+2}E^{3k_{d+1}+1})+
m_3(E^{3k_1}E^{1}E^{3k_2+\cdots+3k_{d+1}+2})
\\&&  +m_2(E^{3k_1}E^{3k_2+\cdots+3k_{d+1}+2}).
\end{eqnarray*}

We complete the proof.
\end{proof}

\begin{lemma}\label{le2}Let $A$ be a bi-Koszul algebra.
Then
$$
\mathbf{ W}^{3k+2} \subseteq \sum_{\substack{i_1+i_2=k,\\
i_1\geq1, i_2\geq0}}m_2(E^{3i_1}E^{3i_2+2})+
\mathbf{U}^{3k+2}+\mathbf{V}^{3k+2}.
$$
\end{lemma}

\begin{proof} By ignoring the lower grading, we may write
$$W^{3k+2}=\sum
m_4(E^{3k_1+i_1}E^{3k_2+i_2}E^{3k_3+i_3}E^{3k_4+i_4})$$ where the
sum runs over all $i_1+i_2+i_3+i_4=7$ ($i_j=1$ or $2$) and $k_j\ge
0$.

First, assume $k_1+k_2+k_3+k_4=0$. In this case, the first or last
component, say the last component, must be $E^2_d$. Using
\textsf{SI(d+3)}, we have
\begin{eqnarray*}
&& m_4(E^{i_1}E^{i_2} E^{i_3}E^{2}_{d})\\
&=&m_4(E^{i_1}E^{i_2}E^{i_3}m_d(E^{1}\cdots E^{1}))
\\&\subseteq& m_{d+1}(E^{4}E^{1}\cdots
E^{1})+m_{d+1}(E^{1}E^{4}\cdots E^{1})+\mathbf{U}^{3k+2}.
\end{eqnarray*}

Next, consider $k_1+k_2+k_3+k_4\geq1$.

If $k_2\geq 1$,
\begin{eqnarray*}
&& m_4(E^{3k_1+i_1}E^{3k_2+i_2}E^{3k_3+i_3}E^{3k_4+i_4})
\\&\subseteq&m_4(E^{3k_1+3k_2+i_1}E^{i_2}
E^{3k_3+i_3}E^{3k_4+i_4})+ \mathbf{U}^{3k+2}+\mathbf{V}^{3k+2}.
\end{eqnarray*}

If $k_3\geq 1$, by the similar method we get
\begin{eqnarray*}
&&m_4(E^{3k_1+i_1}E^{3k_2+i_2}E^{3k_3+i_3}E^{3k_4+i_4})
\\&\subseteq&m_4(E^{3k_1+3k_2+3k_3+i_1}E^{i_2}E^{i_3}E^{3k_4+i_4})
+\mathbf{U}^{3k+2}+ \mathbf{V}^{3k+2}.
\end{eqnarray*}

Now we assume that $k_1\geq1$ (the case of $k_4\geq1$ is
symmetrical). Whether $E^{3k_1+i_1}=E^{3k_1+1}$ or $E^{3k_1+2}$, by
\textsf{SI(5)}  we get
$$
m_4(E^{3k_1+i_1}E^{3k_2+i_2}E^{3k_3+i_3}E^{3k_4+i_4})\subseteq
\mathbf{U}^{3k+2}+\mathbf{V}^{3k+2}.
$$

We complete the proof.\end{proof}

Examining the table in Proposition \ref{pp0} again, the following
result is clear from Proposition \ref{pp2},  the lemmas \ref{le1}
and \ref{le2}.

\begin{corollary}\label{pp03}Let $A$ be a bi-Koszul algebra. Then the
actions of $m_4, m_d, m_{d+1}$ which  fall into $E^{\geq4}$ are
determined by the actions of $m_2$ and $m_3$.     \hfill{$\square$}
\end{corollary}

Now we can state an equivalent statement of the bi-Koszul algebra.

\begin{theorem}\label{thm} Let $A$ be a locally finite, connected
graded algebra generated in degree 1. Then $A$ is a bi-Koszul
algebra if and only if there exists a reduced $A_\infty$-algebra
$(E(A); \{m_i\})$ which is $[m_2, m_3]$-finitely generated by
$E^{1}(A)$, $E^{2}(A)$ and $E^{3}(A)$ with $E^1(A)= E^1_1(A)$,
$E^2(A)= E^2_d(A)\oplus E^2_{d+1}(A)$, $E^3(A)= E^3_{2d}(A)$.
\end{theorem}

\begin{proof}
Assume that  $A$ is a bi-Koszul algebra. Take an $A_\infty$-algebra
$(E(A); \{m_i\})$ that is generated by $E^1(A)$, then $(E(A);
\{m_i\})$ is a reduced $A_\infty$-algebra by Corollary \ref{cc1}. By
Proposition \ref{pp2}, the lemmas \ref{le1} and \ref{le2}, we get
$E(A)$ begins with $E^1(A)= E^1_1(A)$, $E^2(A)= E^2_d(A)\oplus
E^2_{d+1}(A)$, $E^3(A)= E^3_{2d}(A)$, and for each $k\geq 1$,
{\small\begin{eqnarray*} && E^{3k+3}=\sum_{k_1+k_2=k}
m_2(E^{3k_1}E^{3k_2}),\\
&&E^{3k+1}=\sum_{k_1+k_2=k} m_2(E^{3k_1}E^{3k_2+1})
=\sum_{k_1+k_2=k} m_2(E^{3k_1+1}
E^{3k_2}),\\
&&E^{3k+2}=\sum_{k_1+k_2=k} m_2(E^{3k_1}E^{3k_2+2})+\sum_{k_1+k_2=k}
m_2(E^{3k_1+2}E^{3k_2})+\mathbf{U}^{3k+2}+\mathbf{V}^{3k+2},
\end{eqnarray*}}\noindent which implies that $(E(A); \{m_i\})$ is $[m_2,
m_3]$-finitely generated by $E^1(A), E^2(A)$ and $E^3(A)$.

The converse is straightforward by comparing the lower grading.
\end{proof}

The result above tells that the obstruction in Theorem \ref{thm1}
can be described by the  multiplications $m_2$ and $m_3$. Using the
multiplications $m_2$ and $m_3$, we may also give a criteria for a
bi-Koszul algebra to be strongly.

\begin{prop}\label{pp1}Let $A$ be a bi-Koszul algebra. Then $A$  is
strongly if and only if for any $A_\infty$-algebra $(E(A);
\{m_i\})$,
\begin{eqnarray*} && \mathbf{U}^{3k+2}\subseteq
m_2(E^{2}_{d+1}E^{3k}),\;\; and
\\&& \mathbf{V}^{3k+2}\subseteq \sum_{k_1+k_2=k} m_2(E^{3k_1}
E^{3k_2+2})+\sum_{k_1+k_2=k} m_2(E^{3k_1+2}E^{3k_2}).\hskip20mm
\end{eqnarray*}
\end{prop}

\begin{proof} The necessity is obvious. To show the condition being
sufficient, we need only to check
$$E^{3k+2}_{2dk+d+1}=m_2(E^{2}_{d+1}E^{3k}),\quad \mbox{for any}\;\;
k\geq1.\eqno(*)$$ The reason is that the obstruction arises only
from the bigger degree $2dk+d+1$ in $\varDelta_d(3k+2)$ as we
pointed before.


We first show that
$$m_2(E^{2}_{d+1}E^{3k})+\mathbf{U}^{3k+2}=m_2(E^{3k}E^{2}_{d+1})
+\mathbf{U}^{3k+2}.\eqno(**)$$ In fact, by the Stasheff's identity
\textsf{SI(d+2)}
\begin{eqnarray*}
m_2(E^{3k}E^{2}_{d+1})&=&m_2(E^{3k}m_{d+1}(E^{1}\cdots E^{1}))
\\&\subseteq& m_{d+1}(m_2(E^{3k}E^{1})\cdots E^{1})+\mathbf{U}^{3k+2}
\\&=&m_{d+1}(m_2(E^{1}E^{3k})\cdots E^{1})+\mathbf{U}^{3k+2}
\\&\subseteq& m_{d+1}(E^{1}m_2(E^{3k}E^{1})\cdots E^{1})+\mathbf{U}^{3k+2}
\\&\subseteq &\cdots\cdots
\\&\subseteq& m_2(E^{2}_{d+1}E^{3k})+\mathbf{U}^{3k+2},
\end{eqnarray*}
which implies one inclusion relation, the opposite inclusion is
similar to prove.

Using ($**$) and the condition on $\mathbf{U}^{3k+2}$, we have
$m_2(E^{3k}E^{2}_{d+1})\subseteq m_2(E^{2}_{d+1}E^{3k})$. Again from
the conditions on $\mathbf{U}^{3k+2}$, $\mathbf{V}^{3k+2}$ and the
proof of Theorem \ref{thm}, we get
$$
E^{3k+2}=\sum m_2(E^{3k_1}E^{3k_2+2})+\sum m_2(E^{3k_1+2}E^{3k_2}).
$$

We prove ($*$) by induction on $k\ge 1$.
$$
E^{5}_{3d+1}=m_2(E^{3}E^{2}_{d+1})+m_2(E^{2}_{d+1}E^{3})
=m_2(E^{2}_{d+1}E^{3}).
$$
Suppose that $E^{3i+2}_{2di+d+1}=m_2(E^{2}_{d+1}E^{3i})$ for all
$1\leq i< k$. Now, for any $k_1+k_2=k\geq2$ ($k_1\geq1$ and
$k_2\geq1$),
$$
m_2(E^{3k_1}E^{3k_2+2}_{2dk_2+d+1})=m_2(E^{3k_1}m_2(E^{2}_{d+1}E^{3k_2}))
\subseteq m_2(E^{3k_1+2}_{2dk_1+d+1}E^{3k_2})
$$
and
$$
m_2(E^{3k_1+2}_{2dk_1+d+1}E^{3k_2})=m_2(m_2(E^{2}_{d+1}E^{3k_1})E^{3k_2})
\subseteq m_2(E^{2}_{d+1}E^{3k}).
$$
This proves, for any $k\geq1$, $E^{3k+2}_{2dk+d+1}\subseteq
m_2(E^{2}_{d+1}E^{3k})$, the opposite inclusion is trivial, so we
have ($*$).

Hence the bi-Koszul algebra $A$ is strongly.
\end{proof}

\begin{corollary}\label{cor1}
If $A$ is a bi-Koszul algebra with $gl.dim(A)\leq4$, or its
Ext-algebra $E(A)$ with $m_3=0$, then $A$ is a strongly bi-Koszul
algebra.
\end{corollary}

\begin{proof}
It is clear since each assumption implies
$\mathbf{U}^{3k+2}=\mathbf{V}^{3k+2}=0$ for any $k>0$ by Proposition
\ref{pp1}.
\end{proof}
\begin{exa} Any truncated bi-Koszul algebra
is strongly.
\end{exa}

\subsection{Generated by $\mathbf{E^1(A)}$}

Let $A$ be a bi-Koszul algebra, Theorem 3.6 tells that any
$A_\infty$-algebra $E(A)$ is $[m_2, m_3]$-finitely generated by
$E^{1}(A)$, $E^{2}(A)$ and $E^{3}(A)$. On the other hand, Keller has
claimed that there exists an $A_\infty$-algebra structure on $E(A)$
which is generated by $E^1(A)$. In this subsection, we discuss the
universality of the property that $E(A)$ is generated by $E^1(A)$ as
an $A_\infty$-algebra.

For example, let $A$ be an Artin-Schelter regular algebra listed in
\cite[Theorem A]{LPWZ2}, then any $A_\infty$-algebra $E(A)$ is
generated by $E^1(A)$.

Before discussing,  we give a general result which points out that a
strict isomorphism of $A_\infty$-algebras can be obtained from a
quasi-isomorphism of $A_\infty$-algebras. This was found in
\cite{HL} for single $A_\infty$-algebras.

\begin{lemma}\label{le3}
Let $(E; \{m_i\})$ and $(E'; \{m'_i\})$ be two minimal
$A_\infty$-algebras, and $\{f_i\}:(E; \{m_i\})\to (E'; \{m'_i\})$ a
quasi-isomorphism between them. Then
\begin{enumerate}
\item $(E'; \{m''_i\})$ is a minimal $A_\infty$-algebra where
$m''_i:=f_1m_i(f^{-1}_1\otimes\cdots\otimes f^{-1}_1)$ with
$m''_2=m'_2$.
\item $\{g_i\}:(E; \{m_i\})\to (E'; \{m''_i\})$ is a strict isomorphism
of $A_\infty$-algebras  where $g_1=f_1$ and $g_i=0$ for all
$i\geq2$.
\end{enumerate}
\end{lemma}

\begin{proof} Since $m_1=m'_1=0$, $f_1: (E; m_2)\to (E'; m'_2)$ is an
isomorphism. To prove the first statement, we need the Stasheff's
morphism identities \textsf{SI(n)} ($n=1, 2, \cdots$) for
$\{m''_i\}$. Note that the degrees of both $f_1$ and $f^{-1}_1$ are
zero, the Koszul sign convention can be neglected in the following.
For any $i+t+j=n$ and $l=i+1+j$,
\begin{eqnarray*}
&&m''_l(1^{\otimes i}\otimes m''_t \otimes 1^{\otimes j})
\\&=& f_1m_l(f^{-1}_1\otimes\cdots\otimes f^{-1}_1)(1^{\otimes i}\otimes
f_1m_t(f^{-1}_1\otimes\cdots\otimes f^{-1}_1) \otimes 1^{\otimes j})
\\&=&f_1m_l(f^{-1}_1\otimes\cdots\otimes f^{-1}_1)
(1^{\otimes i}\otimes f_1m_t \otimes 1^{\otimes j})(1^{\otimes
i}\otimes f^{-1}_1\otimes\cdots\otimes f^{-1}_1\otimes1^{\otimes j})
\\&=&f_1m_l(1^{\otimes i}\otimes m_t \otimes 1^{\otimes j})
(f^{-1}_1\otimes\cdots\otimes f^{-1}_1),
\end{eqnarray*}
hence
\begin{eqnarray*}
&&\sum_{n=i+t+j}(-1)^{i+tj}m''_l(1^{\otimes i}\otimes m''_t \otimes
1^{\otimes j})
\\&=&\sum_{n=i+t+j}(-1)^{i+tj}f_1(m_l(1^{\otimes i}\otimes m_t
\otimes 1^{\otimes j})) (f^{-1}_1\otimes\cdots\otimes f^{-1}_1)\\
&=&f_1\Big(\sum_{n=i+t+j}(-1)^{i+tj}(m_l(1^{\otimes i}\otimes m_t
\otimes 1^{\otimes j}))\Big) (f^{-1}_1\otimes\cdots\otimes f^{-1}_1)
\\&=&0.
\end{eqnarray*}
Moreover, $m''_2=f_1m_2(f^{-1}_1\otimes f^{-1}_1)=m'_2(f_1\otimes
f_1)(f^{-1}_1\otimes f^{-1}_1)=m'_2$.

Clearly, $f_1m_i=m''_i(f_1\otimes\cdots\otimes f_1)$. So $\{g_i\}$
is a strict isomorphism between $(E; \{m_i\})$ and $(E';
\{m''_i\})$.
\end{proof}

Let $A$ be a bi-Koszul algebra determined by $\varDelta_d$,
$E:=E(A)$ the Ext-algebra of $A$. There is a quasi-isomorphism
between two $A_\infty$-algebra structures on $E(A)$:
$$
\{f_i\}:(E(A); \{m_i\})\rightarrow (E(A); \{m'_i\}).
$$
Now assume that $(E(A); \{m_i\})$ is generated by $E^1(A)$, it is a
natural question that whether the same claim is true for $(E(A);
\{m'_i\})$.

The following facts are immediately:
\begin{itemize}{
\item [(i)] if $n\geq2$, $f_n(E^{1}\cdots E^{1})=0$;
\item [(ii)] if $d\geq3$, $m_2(E^{1}E^{1})=m_2(E^{1}E^{2})
=m_2(E^{2}E^{1})=0$;
\item [(iii)] there are only $m'_d$ and $m'_{d+1}$ whose actions can
fall into $E^{3}$ in $(E(A); \{m'_i\})$;
\item [(iv)] $f_1: E(A)\rightarrow E(A)$ is an isomorphism. }
\end{itemize}

Combining with Lemma \ref{le0} and Proposition \ref{pp0}, we can
write the Stasheff's morphism identities, for small $n$ or in some
special cases, more clearly.
\begin{enumerate}
\item \textsf{MI(2)}:\; $ f_1m_2=m'_2(f_1\otimes f_1)$;
\item \textsf{MI(3)}:\; $f_1m_3+f_2(m_2\otimes1)-f_2(1\otimes m_2)$
\item [] \hskip25mm
$=m'_3(f_1\otimes f_1\otimes f_1)+m'_2(f_1\otimes
f_2)-m'_2(f_2\otimes f_1)$;
\item \textsf{MI(d)}\; acting on $E^{1} \cdots E^{1}$ or $E^{1}\cdots
E^{2}_{d+1}\cdots E^{1}$ can be reduced as
$$f_1m_d=m'_d(f_1\otimes \cdots \otimes f_1);$$
\item \textsf{MI(d+1)}\; acting on $E^{1} \cdots E^{1}$ can be reduced as
$$\qquad\quad f_1m_{d+1}+(-1)^{d}f_2(m_d\otimes1)+(-1)^{d+1}f_2(1\otimes
m_d)=m'_{d+1}(f_1\otimes \cdots \otimes f_1);$$
\item \textsf{MI(d+1)}\; acting on $E^{1} \cdots E^{2}_{d}\cdots E^{1}$
can be reduced as
$$\begin{array}{l}\quad f_1m_{d+1}+(-1)^{d}f_2(m_d\otimes1)+(-1)^{d+1}
f_2(1\otimes m_d)\\=m'_{d+1}(f_1\otimes \cdots \otimes
f_1)+\sum_{1\leq j \leq d}(-1)^{d-j}m'_d(\underbrace{f_1\otimes
\cdots \otimes f_2}_j\otimes\cdots\otimes f_1).\end{array}$$
\end{enumerate}

\begin{prop}\label{pp3}
Assume $E^{3}(A)$ is generated by $E^{1}(A)$ and $E^{2}(A)$ with
$\{m_i\}$. If $f_2(m_d\otimes 1)=f_2(1\otimes m_d)$, then $E^{3}(A)$
is also generated by $E^{1}(A)$ and $E^{2}(A)$ with $\{m'_i\}$.
\end{prop}

\begin{proof}
The hypothesis on $E^{3}(A)$ tells us that
$$
E^{3}=\sum m_d(E^{1}\cdots E^{2}_{d+1}\cdots E^{1})+\sum
m_{d+1}(E^{1}\cdots E^{2}_d\cdots E^{1}).
$$
So  $E^{3}=f_1(E^{3})=\sum f_1m_d(E^{1}\cdots E^{2}_{d+1}\cdots
E^{1})+\sum f_1m_{d+1}(E^{1}\cdots E^{2}_d\cdots E^{1}).$

By the Stasheff's morphism identities listed above, we have
\begin{eqnarray*}
&&f_1m_d(E^{1}\cdots E^{2}_{d+1}\cdots E^{1})
\\&\subseteq& m'_d(f_1\otimes \cdots \otimes f_1)(E^{1}\otimes\cdots\otimes
E^{2}_{d+1}\otimes\cdots \otimes E^{1})
\\&=&m'_d(E^{1}\cdots E^{2}_{d+1}\cdots E^{1}),
\end{eqnarray*}
and
\begin{eqnarray*}
&&f_1m_{d+1}(E^{1}\cdots E^{2}_d\cdots E^{1})
\\&\subseteq&m'_{d+1}(f_1\otimes \cdots \otimes f_1)(E^{1}\otimes\cdots
\otimes E^{2}_d\otimes\cdots \otimes E^{1})
\\& & +\sum_{1\leq j \leq d}(-1)^{d-j}m'_d(f_1\otimes \cdots \otimes f_2
\otimes\cdots\otimes f_1)(E^{1}\otimes\cdots \otimes
E^{2}_d\otimes\cdots \otimes E^{1}).
\\&\subseteq&m'_{d+1}(E^{1}\cdots E^{2}_d\cdots E^{1})
+\sum m'_d(E^{1}\cdots E^{2}_{d+1}\cdots E^{1}).
\end{eqnarray*}
Thus, $E^{3}\subseteq \sum m'_d(E^{1}\cdots E^{2}_{d+1}\cdots
E^{1})+\sum m'_{d+1}(E^{1}\cdots E^{2}_d\cdots E^{1})$. We complete
the proof.
\end{proof}

Now, we can prove the main theorem of this subsection.

\begin{theorem} \label{thm2} Let $\{f_i\}:(E(A); \{m_i\})\to
(E(A); \{m'_i\})$ be a quasi-isomorphism with $f_2(m_i\otimes
1)=f_2(1\otimes m_i)$ for $i=2, d$, suppose that $(E(A); \{m_i\})$
is generated by $E^1(A)$. Then $(E(A); \{m'_i\})$ is also generated
by $E^1(A)$.
\end{theorem}

\begin{proof}Since $m_1=m'_1=0$, $f_1:(E(A), m_2)\rightarrow (E(A), m'_2)$
is an isomorphism with degree zero.  By the proof of Theorem
\ref{thm} and \textsf{MI(2)}, we have
\begin{eqnarray*}
&E^{3k+3}&= f_1(E^{3k+3})=f_1(\sum m_2(E^{3k_1}E^{3k_2}))
\\&&\subseteq  \sum m'_2(f_1\otimes f_1)(E^{3k_1}\otimes E^{3k_2})=
\sum m'_2(E^{3k_1}E^{3k_2});\\
&E^{3k+1}&=f_1(E^{3k+1})=f_1(\sum m_2(E^{3k_1+1}E^{3k_2}))
\\&&\subseteq  \sum m'_2(f_1\otimes f_1)(E^{3k_1+1}\otimes E^{3k_2})=
\sum m'_2(E^{3k_1+1}E^{3k_2});\\
&E^{3k+2}&=f_1(E^{3k+2})\\&&=f_1(\sum m_2(E^{3k_1}E^{3k_2+2})+\sum
m_2(E^{3k_1+2}E^{3k_2})+\mathbf{U}^{3k+2}+\mathbf{V}^{3k+2})
\\&&\subseteq \sum m'_2(E^{3k_1}E^{3k_2+2})+\sum
m'_2(E^{3k_1+2}E^{3k_2})+f_1(\mathbf{U}^{3k+2})+f_1(\mathbf{V}^{3k+2}).
\end{eqnarray*}

For any $E^{3i_1}E^{3i_2+1}E^{3i_3+2}_{2di_3+d}$ ($i_1+i_2+i_3=k,
i_1\geq1,i_2,i_3\geq0$), the assumption
$f_2(m_2\otimes1)=f_2(1\otimes m_2)$ and \textsf{MI(3)} imply that
\begin{eqnarray*}
&&f_1m_3(E^{3i_1}E^{3i_2+1}E^{3i_3+2}_{2di_3+d})
\\&\subseteq& m'_3(E^{3i_1}E^{3i_2+1}E^{3i_3+2}_{2di_3+d})+
m'_2(f_2(E^{3i_1}E^{3i_2+1})f_1(E^{3i_3+2}_{2di_3+d}))
\\&&+m'_2(f_1(E^{3i_1})f_2(E^{3i_2+1}E^{3i_3+2}_{2di_3+d}))
\\&\subseteq& m'_3(E^{3i_1}E^{3i_2+1}E^{3i_3+2}_{2di_3+d})+m'_2( E^{3i_1
+3i_2} E^{3i_3+2}_{2di_3+d})
\\&&+ m'_2(E^{3i_1}E^{3i_2+3i_3+2}_{2di_2+2di_3+d}).
\end{eqnarray*}

By the same method, we obtain that the action  of $f_1$ on every
component of $\mathbf{U}^{3k+2}$ or $\mathbf{V}^{3k+2}$ falls into
$\sum m'_2(E^{3k_1}E^{3k_2+2})+\sum
m'_2(E^{3k_1+2}E^{3k_2})+\mathbf{U}'^{3k+2}+\mathbf{V}'^{3k+2}$
where $\mathbf{U}'^{3k+2}$ and $\mathbf{V}'^{3k+2}$  correspond to
$\mathbf{U}^{3k+2}$ and $\mathbf{V}^{3k+2}$, respectively, changing
the multiplication $m_3$ to $m'_3$. Thus,
$$
E^{3k+2}=\sum m'_2(E^{3k_1}E^{3k_2+2})+\sum
m'_2(E^{3k_1+2}E^{3k_2})+\mathbf{U}'^{3k+2}+\mathbf{V}'^{3k+2}.
$$
Since $E^{2}$ is generated by $E^{1}$, and $E^{3}$ is generated by
$E^{1}$ and $E^{2}$ by Proposition \ref{pp3}, we get $(E(A);
\{m'_i\})$ is also generated by $E^1$.
\end{proof}

\begin{corollary}Assume the bi-Koszul algebra $A$ is strongly.  If
$(E(A); \{m_i\})$ is generated by $E^1$, and $f_2(m_d\otimes
1)=f_2(1\otimes m_d)$, then $(E(A); \{m'_i\})$ is also generated by
$E^1$.\qed
\end{corollary}

Continuing to consider the quasi-isomorphism between two
$A_\infty$-structures on $E(A)$, $ \{f_i\}:(E(A);
\{m_i\})\rightarrow (E(A); \{m'_i\}), $ some extra hypothesis will
make $\{f_i\}$ to be a strict isomorphism which guarantees the
properties of such two $A_\infty$-algebras identify with each other.

\begin{theorem}\label{thm3}
If $f_2(m_i\otimes1)=f_2(1\otimes m_i)$ for $i=2, d$, and
$$
\sum_{1\leq j\leq i} (-1)^{i-j}m'_i(\underbrace{f_1\otimes
\cdots\otimes f_2}_j \otimes\cdots f_1)=0,
$$
then $g=f_1:(E(A); \{m_i\})\rightarrow (E(A); \{m'_i\})$ is a strict
isomorphism.
\end{theorem}

\begin{proof}
By Lemma \ref{le3}, $g=f_1:(E(A); \{m_i\})\rightarrow (E(A);
\{m''_i\})$ is a strict isomorphism and $m''_2=m'_2=m_2$.  The
assumption directly implies $m''_3=f_1m_3(f_1^{-1}\otimes
f_1^{-1}\otimes f_1^{-1})=m'_3$. For any $x_1,\cdots,x_d\in E^{1}$,
or one of $x_j$'s in $E^{2}_{d+1}$ and the others in $E^{1}$, we
have
$$m''_d(x_1,\cdots,x_d)=f_1m_d(f_1^{-1}(x_1),\cdots,f_1^{-1}(x_d))
=m'_d(x_1,\cdots,x_d)$$ by the Stasheff's morphism identities listed
above. So $m''_d|_{(E^{1})^{\otimes d}}=m'_d|_{(E^{1})^{\otimes d}}$
and $m''_d|_{E^{1}\cdots E ^{2}_{d+1}\cdots E^{1}
}=m'_d|_{E^{1}\cdots E ^{2}_{d+1}\cdots E^{1} }$. By the same
method, we check that $m''_{d+1}|_{(E^{1})^{\otimes
(d+1)}}=m'_{d+1}|_{(E^{1})^{\otimes (d+1)}}$ and
$m''_{d+1}|_{E^{1}\cdots E ^{2}_d\cdots E^{1}
}=m'_{d+1}|_{E^{1}\cdots E ^{2}_d\cdots E^{1} }$. Thus, $(E(A);
\{m'_i\})$ and $(E(A); \{m''_i\})$ are the same. The result follows
immediately.
\end{proof}

\begin{corollary}Assume the bi-Koszul algebra $A$ is strongly.
If $f_2(m_d\otimes1)=f_2(1\otimes m_d)$ and  $\sum_{1\leq j\leq d}
(-1)^{d-j}m'_d(\underbrace{f_1\otimes \cdots\otimes f_2}_j
\otimes\cdots \otimes f_1)=0$. Then
$$g=f_1:(E(A);
\{m_i\})\rightarrow (E(A); \{m'_i\})$$ is a strict isomorphism.
\end{corollary}

\begin{proof}
By the proof of Theorem \ref{thm3}.
\end{proof}

\vskip1cm
\bibliographystyle{amsplain}

\begin{thebibliography}{11}
\bibitem{BGS}
A. Beilinson, V. Ginzburg and W. Soergel, \emph{Koszul duality
patterns in representation theory}, J. Amer. Math. Soc. \textbf{9}
(1996), no. 2, 473--52
\bibitem{Be} R. Berger, \emph{Koszulity of nonquadratic algebras}, J.
Algebra, \textbf{239} (2001), 705-734.
\bibitem{CS} T. Cassidy and B. Shelton, \emph{Generalizing the notion of a
Koszul algebra}, Math. Z., \textbf{260} (2008),  93-114.
\bibitem{HL} J.-W. He and D.-M. Lu, \emph{Higher Koszul algebras and
$A$-infinity algebras}, J. Algebra, \textbf{293} (2005), 335-362.
\bibitem{Ka}
T. V. Kadeishvili,  \emph{On the theory of homology of fiber
spaces}, (Russian) International Topology Conference (Moscow State
Univ., Moscow, 1979). Uspekhi Mat. Nauk \textbf{35} (1980), no. 3
(213), 183--188. The English translation was published in Russian
Math. Surveys, \textbf{35}, no.3 (1980), 231-238.
\bibitem{K1} B. Keller, \emph{$A$-infinity algebras in representation
theory}, Contribution to the Proceedings of ICRA IX. Beijing: Peking
University Press, 2000.
\bibitem{K2} B. Keller,  \emph{Introduction to $A$-infinity algebras and
modules}, Homology Homotopy Appl., \textbf{3} (2001), 1-35
(electronic).
\bibitem{LPWZ1} D.-M. Lu, J. H. Palmieri, Q.-S. Wu and J. J. Zhang,
\emph{ A-infinity algebras for ring theorists}, Algebra Colloq.,
\textbf{11} (2004), 91-128.
\bibitem{LPWZ2} D.-M. Lu, J. H. Palmieri, Q.-S. Wu and J. J. Zhang,
\emph{Regular algebras of dimension 4 and their
$A_\infty$-Ext-algebras}, Duke Math. J., \textbf{137} (2007),
537-584.
\bibitem{LS} D.-M. Lu and J.-R. Si,  \emph{Koszulity of algebras with non-pure
resolutions}, Comm. Alg., accepted for publication.
\bibitem{LHL} J.-F. L\"{u}, J.-W. He and D.-M. Lu,  \emph{Piecewise-Koszul
algebras}, Sci. China Ser. A, \textbf{50} (2007), 1795-1804.
\bibitem{P} S. Priddy, \emph{Koszul
resolutions}, Trans. Amer. Math. Soc., \textbf{152} (1970), 39-60.
\end{thebibliography}

\end{document}